\documentclass[11pt]{amsart}
\usepackage{graphicx}
\usepackage{amsmath}
\usepackage{amscd}
\usepackage{mathrsfs}
\usepackage{latexsym}
\usepackage{epstopdf}

\newcommand{\Q}{\mathbb{Q}}
\newcommand{\Z}{\mathbb{Z}}
\newcommand{\C}{\mathbb{C}}

\newcommand{\set}[1]{\{ #1 \}}

\newcommand{\dbar}{\overline{\partial}}
\newcommand{\twopi}[1]{\frac{ #1 }{2\pi i}}
\newcommand{\Proj}{\mathbb{P}}
\newcommand{\Bl}[1]{\widetilde{#1}}

\newcommand{\jacob}[2]{\frac{\vartheta(\twopi{#1}-#2 z)\vartheta(z)}{\vartheta(\twopi{#1}-z)\vartheta(#2 z)}}

\newcommand{\ellipstd}[1]{\frac{#1\vartheta(\frac{#1}{2\pi i}-z,\tau)}{\vartheta(\frac{#1}{2\pi i},\tau)}}
\newcommand{\ellnorm}[1]{\frac{\vartheta(#1-z,\tau)}{\vartheta(#1,\tau)}}
\newcommand{\jac}[2]{\frac{\vartheta(#1-(#2+1)z,\tau)\vartheta(z,\tau)}{\vartheta(#1-z,\tau)\vartheta((#2+1)z,\tau)}}
\newcommand{\jacobtau}[2]{\frac{\vartheta(\twopi{#1}-#2 z,\tau)\vartheta(z,\tau)}{\vartheta(\twopi{#1}-z,\tau)\vartheta(#2 z,\tau)}}

\newtheorem{thm}{Theorem}
\newtheorem{prop}{Proposition}
\newtheorem{lem}{Lemma}

\newtheorem{defn}{Definition}
\newtheorem{rmk}{Remark}
\newtheorem{cor}{Corollary}

\title{Singular McKay correspondence for normal surfaces}
\author{Robert Waelder}
\address{rwaelder@math.uic.edu}
\thanks{The author supported by Clay Liftoff Fellowship and NSF Post-doctoral Fellowship}
\begin{document}
\begin{abstract}
We define the singular orbifold elliptic genus and $E$-function for all normal surfaces without strictly log-canonical singularities, and prove the analogue of the McKay correspondence in this setting. Our invariants generalize the stringy invariants defined by Willem Veys for this class of singularities. We show that the ability to define these invariants is closely linked to rigidity phenomena associated to the elliptic genus.
\end{abstract}

\maketitle

\section{Introduction}
Suppose $Y=X/G$ is the quotient of a smooth variety $X$ by a finite group $G$. There are two natural procedures for defining smooth invariants, such as Chern numbers and Hodge numbers, on the singular variety $Y$. First, one could resolve the singularities of $Y$ and define smooth invariants on $Y$ by using the data of the resolution. Alternatively, one could regard $Y$ as an orbifold and define smooth invariants on $Y$ using the orbifold data of $(X,G)$. The observation that both procedures tend to be equivalent is the one interpretation of the {McKay correspondence}.

In this context, we will focus on generalizations of the $E$-polynomial and complex elliptic genus to singular varieties. The former encodes Hodge structure data, the latter encodes Chern number data. Both invariants are natural objects of interest in birational geometry. For example, the $E$-polynomial is a universal object in the Grothendieck ring of varieties. And, by a result of Totaro, the elliptic genus encodes precisely those Chern numbers that are invariant under flops. \cite{Totaro} As seen in the work of Batyrev, Borisov, Libgober, and Chin-Lung Wang \cite{Batyrev,BL_Sing,CLW}, in order to generalize these invariants to the birational category, one has to extend the objects of study to include divisor pairs $(X,D)$ which are Kawamata log-terminal. This just means that $K_X-D$ is $\Q$-Cartier, and that the discrepancy coefficients of any resolution of $(X,D)$ are all greater than $-1$.
Given these technical constraints on the singularities of $(X,D)$, both the $E$-polynomial and elliptic genus of $(X,D)$ turn out to be functorial with respect to birational morphisms of pairs. This functoriality is in fact essential for ensuring that the given invariant is well-defined. 

In this paper we will examine the following singular analogue of the McKay correspondence: Let $X$ be a {singular} variety with a finite group action. As before, there are two natural procedures for defining smooth invariants on the quotient $Y=X/G$: $(1)$ we could resolve the singularities of $Y$ and define smooth invariants on $Y$ using the data of the resolution; $(2)$ we could construct a $G$-equivariant resolution $\Bl{X}\rightarrow X$, and define smooth invariants on $Y$ using the orbifold data of $(\Bl{X},G)$. If $X$ has at worst log-terminal singularities, it is known that these two procedures are equivalent. 

Unfortunately, for non-log-terminal singularities, poles appear in both the Batyrev and Borisov-Libgober-Wang expressions for the $E$-function and elliptic genus. It is therefore unclear whether the above-described program even makes sense for varieties with worse than log-terminal singularities. At least for the non-orbifold case (i.e., procedure $(1)$ described above) Willem Veys \cite{Veys} has succeeded in extending Batyrev's $E$-function to ``almost all" normal surfaces without log-terminal singularities. Moreover this author \cite{RW_Stringy} has extended Borisov and Libgober's singular elliptic genus to the same class of surface singularities investigated by Veys. It is therefore natural to conjecture that the same extensions can be carried out for the orbifold case (procedure $(2)$), and that the McKay correspondence can be verified in this setting. 

In what follows, we will verify this conjecture, and along the way reveal an interesting relationship between the extension of smooth invariants to non-log-terminal singularities and certain rigidity theorems on toric varieties. 

\section{Preliminaries}
\subsection{Singularities from the minimal model program}
Let $Y$ be a normal surface with $\Q$-Cartier canonical divisor, and $f:X\rightarrow Y$ a resolution of singularities. We say that $X$ is a \emph{log resolution} if the exceptional curves of $f:X\rightarrow Y$ are smooth with simple normal crossings. If $\cup_I E_i$ is the exceptional locus of $f$, we define the discrepancies $a_i$ of $E_i$ by the formula:
\begin{align*}
K_X = f^*K_Y + \sum_{I} a_i E_i
\end{align*}
More generally, for \emph{any} normal surface $Y$ with $\Q$-Weil divisor $\Delta \subset Y$, and log resolution $f:X\rightarrow Y$, we may define the pullback $f^*(K_Y-\Delta)$ by the following procedure given by Mumford: Let $\Bl{\Delta}$ denote the proper transform of $\Delta$ with divisor coefficients equal to the corresponding coefficients of $\Delta$. If we were in the smooth setting, we would define the discrepancy coefficients $a_i$ of $E_i$ by the formula:
\begin{align*}
K_X = f^*(K_Y-\Delta)+\Bl{\Delta}+\sum_I a_i E_i
\end{align*}
Multiplying both sides by $E_j$ gives a system of linear equations
\begin{align*}
K_X E_j -\Bl{\Delta}E_j = \sum_I a_i E_i E_j
\end{align*}
Since the matrix of intersection numbers $\set{E_i E_j}$ is negative definite, we can always find a unique solution $\set{a_i}$ over the rationals. Note that the above equations make sense for any normal surface $Y$ and $\Q$-Weil divisor $\Delta$. We therefore define $f^*(K_Y-\Delta) = K_X-\Bl{\Delta}-\sum_{I}a_i E_i$, where the $a_i$'s are given by the same set of equations. 

Thus, let $Y$ be a normal surface with $\Q$-Weil divisor $\Delta$. Let $f:X\rightarrow Y$ be a log resolution with exceptional components $\cup E_i$. Let $D =\sum a_i D_i$ be the $\Q$-Cartier divisor on $X$ supported on $\cup E_i\cup\Bl{\Delta}$, whose coefficients $a_i$ satisfy the equation $K_X-D = f^*(K_Y-\Delta)$. We refer to $f:(X,D)\rightarrow (Y,\Delta)$ as a log resolution of $(Y,\Delta)$.
The singularities of $(Y,\Delta)$ said to be:\\

\noindent $(1)$ \emph{log-terminal} if $a_i > -1$ for all $i$.\\ 
$(2)$ \emph{log-canonical} if $a_i\geq -1$ for all $i$.\\
$(3)$ \emph{strictly log-canonical} if $a_i\geq -1$ and at least some $a_i = -1$. 
\\

\noindent These classifications are independent of the choice of log resolution. 

\begin{rmk}\rm
More generally, for any regular map between normal surfaces $f:X\rightarrow Y$, we will find it convenient to write $f:(X,D)\rightarrow (Y,\Delta)$ if the pull-back of $K_Y-\Delta$ by $f$ ``makes sense", and $K_X-D = f^*(K_Y-\Delta)$ as $\Q$-Weil divisors.
\end{rmk}

\subsection{Stringy invariants}
Let $(X,D)$ be a smooth log-terminal pair, i.e., $D=\sum_I a_i D_i$ is a sum of smooth divisors with simple normal crossings with $a_i>-1$. Batyrev \cite{Batyrev} defines the stringy $E$-function of $(X,D)$ as follows: For every subset $J\subset I$, let $D_J^o =\cap_J D_j\backslash \cup_{J^c}D_i$. By definition, we set $D_{\emptyset}^o = X\backslash \cup_I D_i$. Then
\begin{align*}
E_{str}(X,D;u,v) =\sum_{J\subset I}E(D_J^o;u,v)
\prod_J\frac{uv-1}{(uv)^{a_i+1}-1}
\end{align*}
Here $E(D_J^o;u,v)$ is the usual $E$-polynomial of $D_J^o$ given by the Deligne mixed Hodge structure on $H^*_c(D_J^o,\C)$. 

The above definition of $E_{str}$ may be interpreted as a motivic integral over the space $J_{\infty}(X)$ if germs of arcs on $X$ (see \cite{Batyrev}). The change of variable formula from motivic integration implies that $E_{str}(X,D;u,v)$ is functorial with respect to birational morphisms of divisor pairs. That is, if $f:(X,D)\rightarrow (Y,\Delta)$ is a birational morphism of smooth log-terminal pairs, then $E_{str}(X,D;u,v)=E_{str}(Y,\Delta;u,v)$. This allows us to define the stringy $E$-function for any log-terminal pair $(Y,\Delta)$, so long as $K_Y-\Delta$ is $\Q$-Cartier (in the lingo of the minimal model program, $(Y,\Delta)$ is \emph{Kawamata log-terminal}). We define $E_{str}(Y,\Delta;u,v) = E_{str}(X,D;u,v)$, for $(X,D)\rightarrow (Y,\Delta)$ a log resolution. The definition is well-defined, since for any two log resolutions $(X_1,D_1)$ and $(X_2,D_2)$, we can find a common resolution $(M,D)$ making the following diagram commute:
$$\begin{CD}
(M,D) @>>> (X_2,D_2)\\
@VVV       @VVV\\
(X_1,D_1)       @>>> (Y,\Delta)
\end{CD}$$
We refer to the specializations $E_{str}(Y,\Delta;u,1)$ as the stringy $\chi_y$ genus of $(Y,\Delta)$ and $E_{str}(Y,\Delta;1,1)=e_{str}(Y,\Delta)$ as the stringy Euler number.

In the situation where $Y = X/G$, where $G$ is a finite group acting on $X$, one can approach the problem of defining the $E$-function of $Y$ by making use of the orbifold data associated to $(X,G)$. The resulting invariant is called the orbifold $E$-function of $(X,G)$. For our purposes, it is convenient to define the orbifold $E$-function for triples $(X,D,G)$, which satisfy the following condition, known as $G$-normality:

\begin{defn}\rm
Let $X$ be a smooth variety, and $D=\cup D_i$ a smooth divisor with simple normal crossings. Let $G$ be a finite group acting holomorphically on $X$ and leaving $D$ invariant. Then $(X,D)$ is $G$-normal if for every $p\in X$, and $D_i$ containing $p$, $\mathrm{stab}_G(p)D_i = D_i$.
\end{defn} 

Let $(X,D=\sum_I a_i D_i)$ be a smooth $G$-normal pair. Following Batyrev, we proceed to define the orbifold $E$-function $E_{orb}(X,D,G;u,v)$. Fix $g\in G$, and let $X^g\subset X$ be a fixed component of $g$. Let $K\subset I$ index the divisors $D_k$ containing $X^g$. Then the normal bundle to $X^g$ splits into a direct sum of character sub-bundles $\bigoplus N_\alpha\oplus\bigoplus_{K}\mathcal{O}(D_k)|_{X^g}$, where $g$ acts on $N_\alpha$ by $e^{2\pi i\alpha}$ and acts on $\mathcal{O}(D_k)|_{X^g}$ by $e^{2\pi i\alpha_k}$. Define the fermionic shift of $g$ as
\begin{align*}
F(g,D) = \sum_\alpha \alpha\cdot\mathrm{rk}(N_\alpha)+\sum_k (1+a_k)\alpha_k
\end{align*}
For $D=0$ this is the usual expression for the Fermionic shift which comes from physics. Let $I^g\subset I$ index the set of $g$-invariant components of $D$. For $J\subset I^g$, let $C(g,J,X^g)$ denote the subgroup of the centralizer of $g$ that leaves $X^g$ and $J$ invariant. Then $E_{orb}(X,D,G;u,v)$ is defined as
\begin{align*}
\sum_{[g],[X^g]}(uv)^{F(g,D)}\sum_{J\subset I^g}
E(X^g\cap D_J^o/C(g,J,X^g);u,v)\prod_J\frac{uv-1}{(uv)^{a_j+1}-1}
\end{align*}
The first two sums run over the conjugacy classes $[g]$ of $G$ and representatives $[X^g]$ of $C(g)$ orbits of components of $X^g$. As for the stringy $E$-function, Batyrev has shown that $E_{orb}(X,D,G;u,v)$ is functorial with respect to $G$-equivariant birational morphisms of $G$-normal pairs.

The relationship between $E_{str}(X/G;u,v)$ and $E_{orb}(X,G;u,v)$ is the subject of the McKay correspondence, which we will discuss at the end of this section. 

Just as the stringy $E$-function provides a means of keeping track of Hodge number data under birational morphisms, a similar procedure exists for studying Chern numbers, which we now review.
 
\subsection{Elliptic genera}

The elliptic genus $Ell(X;z,\tau)$ of an almost complex manifold $X$ is defined as the index of the operator:
$$\dbar\otimes\bigotimes_{n=1}^{\infty}\Lambda_{-yq^{n-1}}T''X\otimes
\Lambda_{-y^{-1}q^n}T'X\otimes S_{q^{n}}T''X\otimes S_{q^n}T'X.$$
Here $T'X$ and $T''X$ are the holomorphic and anti-holomorphic tangent bundles; $\Lambda_t(E)$ and $S_t(E)$ denote the formal sums of exterior and symmetric powers of $tE$; $y = e^{2\pi iz}$ and $q = e^{2\pi i\tau}$. By Riemann-Roch, this index is given by the following integral involving the formal Chern roots $x_i$ of $T'X$:
$$\int_X \prod_{T'X}\frac{x_i\vartheta(\twopi{x_i}-z,\tau)}{\vartheta(\twopi{x_i},\tau)}$$
Here $\vartheta(t,\tau)$ is the Jacobi theta function.

One sees from the above expression that special values of the elliptic genus produce many interesting geometric invariants. For instance, $Ell(X;z,q=0) = y^{-n/2}\chi_{y}(X)$, where $\chi_{y}$ is the Hirzebruch $\chi_{y}$ genus, and $Ell(X;\frac{1}{2},\tau)$ reproduces the signature.

For $(X,D=\sum a_i D_i)$ a smooth log-terminal pair, Borisov, Libgober, and Chin-Lung Wang \cite{BL_Sing, CLW} have defined the elliptic genus of $(X,D)$ by the formula:
\begin{align*}
Ell(X,D;z,\tau)=\int_X \prod_{TX}\frac{\twopi{x_j}\vartheta(\twopi{x_j}-z,\tau)}
{\vartheta(\twopi{x_j},\tau)}\prod_i\jacobtau{D_i}{(a_i+1)}
\end{align*}
In the above formula, $D_i$ represent the classes $c_1(\mathcal{O}(D_i))$. As in the case of the $E$-function, $Ell(X,D;z,\tau)$ is functorial with respect to birational morphisms of pairs. If $(Y,\Delta)$ is Kawamata log-terminal, we may therefore define the singular elliptic genus $Ell(Y,\Delta;z,\tau)$ to be the elliptic genus $Ell(X,D;z,\tau)$ for $(X,D)\rightarrow (Y,\Delta)$ a log-resolution. By the weak factorization theorem \cite{W}, any two log-resolutions of $(Y,\Delta)$ factor into a sequence of blow-ups and blow-downs. Thus, functoriality of $Ell(X,D;z,\tau)$ with respect to birational morphisms ensures that this definition is well-defined.

For $(X,D)$ a $G$-normal pair, following \cite{BL_McKay}, one can define the orbifold elliptic genus $Ell_{orb}(X,D,G;z,\tau)$ in a manner analogous to the definition given by Batyrev for the orbifold $E$-function. Though we will not need to do explicit calculations involving the orbifold elliptic genus, we include its definition here for completeness:

For $g,h \in G$ a 
  commuting pair, let $\set{X^{g,h}_\gamma}$ denote the connected components of their
  common fixed point locus. Fix one such component $X^{g,h}_\gamma$.
  The normal bundle $N_{X^{g,h}_\gamma}$ splits as a sum $\oplus_{\alpha} 
  N_{\alpha}$ over irreducible characters for the subgroup $(g,h)$. 
  For $x \in (g,h)$, let $\alpha(x) \in \Q\cap [0,1)$ be the rational 
  number such that $x$ acts on the fibers of $N_{\alpha}$ as 
  multiplication by $e^{2\pi i\alpha(x)}$. 
  Fix also an irreducible component $D_i$ of $D$. If $X^{g,h}_\gamma\subset 
  D_i$ then $x\in (g,h)$ acts on the fibers of
  $\mathcal{O}(D_i)|_{X^{g,h}_\gamma}$ as multiplication by $e^{2\pi 
  i\epsilon_i(x)}$ for some rational number $\epsilon_i(x)\in \Q\cap 
  [0,1)$. If $X^{g,h}_\gamma$ is not contained in $D_i$, we define 
  $\epsilon_i = 0$. Of course the functions $\alpha$ and $\epsilon_i$ 
  depend on the choice of the commuting pair $(g,h)$ and on the
  connected component $X^{g,h}_\gamma$ of $X^{g,h}$. We will 
  omit making explicit reference to this dependence in order to simplify the 
  notation. The orbifold elliptic genus of $(X,D,G)$ is then given by the formula:
  $$\frac{1}{|G|}\sum_{gh=hg,\gamma}\int_{X^{g,h}_\gamma}
  \prod_{TX^{g,h}_\gamma}\ellipstd{x_j}\times$$
  $$\prod_{N_{\alpha}}\ellnorm{\frac{x_\alpha}{2\pi 
  i}+\alpha(g)-\alpha(h)\tau}e^{2\pi i\alpha(h)z}\times$$
  $$\prod_{D_i}\jac{\frac{D_i}{2\pi 
  i}+\epsilon_i(g)-\epsilon_i(h)\tau}{a_i}e^{-2\pi 
  i\delta_i\epsilon_i(h)z}
  $$
  Here $x_j$ are the Chern roots of $TX^{g,h}_\gamma$, $x_\alpha$ are the Chern roots of $N_\alpha$, and $D_i = c_1(\mathcal{O}(D_i))$.
  One can show that the above formula coincides with Batyrev's orbifold $\chi_y$ genus in the limit $\tau\to i\infty$ by making use of the Lefschetz fixed point formula. 
  
  Finally, if $(X,D)$ admits a torus action $T$ that commutes with $G$, we define the equivariant orbifold elliptic genus $Ell_{orb}(X,D,G;\vec{t},z,\tau)$ by substituting all Chern data in the above formula by their equivariant extensions in the equivariant cohomology ring $H^*_T(X)$. See \cite{RW_McKay} or \cite{LLZ} for details.

As proven in \cite{BL_McKay} (see \cite{RW_McKay} for the $T$-equivariant case), the orbifold elliptic genus is functorial with respect to equivariant birational morphisms of $G$-normal pairs. Since any two $G$-normal resolutions of a $G$-Kawamata log-terminal pair $(X,D)$ can be connected by a sequence of $G$-normal blow-ups and blow-downs (see \cite{BL_McKay}), this functoriality allows us to define the \emph{singular orbifold elliptic genus} of $(X,D,G)$ by following the same procedure as above for the singular elliptic genus. As remarked above, the orbifold $E$-function is also functorial with respect to equivariant birational morphisms of $G$-normal pairs. We therefore define the \emph{singular orbifold $E$-function} similarly.

\subsection{The McKay correspondence}\label{def McKay}
We first review the McKay correspondence in the log-terminal setting for the $E$-function and the elliptic genus. Let $X$ be a $\Q$-Gorenstein normal variety with log-terminal singularities, and let $G$ be a finite group acting on $X$. Let $g:X\rightarrow X/G$ be the quotient map. By the ramification formula, there exists a unique $\Q$-Weil divisor $\Delta_{X/G}\subset X/G$ satisfying $K_X = g^*(K_{X/G}-\Delta_{X/G})$.
As $X$ and $X/G$ are possibly singular, we should clarify what we mean by the pullback of a $\Q$-Weil divisor on $X/G$. Since $X$ is normal, $X/G$ is normal. In particular, $X/G$ is smooth away from a codimension $2$ subset. Let $U$ denote the smooth locus of $X/G$. Since $g:X\rightarrow X/G$ is a finite morphism of normal varieties, $X\backslash g^{-1}(U)$ has codimension at least $2$ in $X$. Thus, if $D$ is any $\Q$-Weil divisor on $X/G$, we may define $g^*D$ by first restricting it to $U$, pulling it back by $g$, and then extending it uniquely to all of $X$. Note, moreover, that the normality of $X$ guarantees that $(X/G,\Delta_{X/G})$ is $\Q$-Cartier. Given this clarification we have the following singular version of the McKay correspondence, which is implicit in the proofs of Batyrev \cite{Batyrev}, and Borisov-Libgober \cite{BL_McKay}:
\begin{thm}
\begin{align*}
&E_{orb}(X,G;u,v) = E_{str}(X/G,\Delta_{X/G};u,v)\\
&Ell_{orb}(X,G;z,\tau) = Ell(X/G,\Delta_{X/G};z,\tau)
\end{align*}
\end{thm}
Note that the McKay correspondence for the $T$-equivariant orbifold elliptic genus is proven in \cite{RW_McKay}.

Our goal in this paper is to prove the above theorem for surfaces without strictly log-canonical singularities. Since we can use the Mumford trick to pull back any $\Q$-Weil divisor on a normal surface by a resolution of singularities, we do not require the $\Q$-Gorenstein condition.

\subsection{Generalization to non-log-terminal singularities}
Notice that we would have to divide by zero in the expressions for the $E$-function and elliptic genus of a $(X,D)$ if any of the divisor coefficients $a_i = -1$. This is the reason for the log-terminality condition in the definitions for the stringy $E$-function and singular elliptic genus. In fact, even if $(Y,\Delta)$ possessed a log-resolution $(X,D)$ with no $-1$ discrepancies, though with some discrepancy coefficients $a_i < -1$, one could still not conclude that $E_{str}(Y,\Delta;u,v)$ and $Ell(Y,\Delta;z,\tau)$ were well-defined. For if $(X',D')$ were another log-resolution with no $-1$ discrepancies, it may happen that a $-1$ discrepancy appeared somewhere along the chain of blow-ups and blow-downs connecting $(X,D)$ to $(X',D')$. Consequently we would have no way to compare, for example, $E_{str}(X,D;u,v)$ with $E_{str}(X',D';u,v)$.

Suppose, however, that $Y$ is a normal surface without strictly log-canonical singularities and $(X,D)\rightarrow (Y,0)$ is the minimal log-resolution. Then Willem Veys \cite{Veys} observed that the only components of $D$ with $-1$ discrepancies are smooth rational curves that intersect at most two other curves at a single point. In general if $D_t\cong \Proj^1$ has coefficient $a_t\neq -1$ and intersects at most two other divisor components $D_{t_k}$ at a single point, then the contribution of $D_t$ to $E_{str}(X,D;u,v)$ is easily computed as:
\begin{align*}
\frac{(uv-1)^2((uv)^{(m_t-1)(a_t+1)}+\ldots + 1)}
{((uv)^{a_{t_1}+1}-1)((uv)^{a_{t_2}+1}-1)}
\end{align*}
where $m_t = -D_t\cdot D_t$. Since the above expression makes sense even in the limit as $a_t\to -1$, Veys defined the contribution to $E_{str}(X,D;u,v)$ coming from each $-1$ discrepancy to be equal to  
\begin{align*}
\frac{m_t(uv-1)^2}
{((uv)^{a_{t_1}+1}-1)((uv)^{a_{t_2}+1}-1)}
\end{align*}
This definition of $E_{str}(X,D;u,v)$ turns out to be functorial with respect to blow-ups, so long as the blow-ups do not occur at a generic point of $D_t$.

It turns out that we can associate to $(X,D)$ a canonical $\Q$-cone of divisors $\set{\Delta}$ so that $E_{str}(X,D;u,v)$ is equal to $\lim_{\varepsilon\to 0}E_{str}(X,D+\varepsilon\Delta';u,v)$ for any $\Delta'\in\set{\Delta}$. We will discuss this issue in detail in section (\ref{orb inv}), and show how we can use this observation to define the orbifold $E$-function and orbifold elliptic genus in this setting. The idea is based on the observation that every component $D_t$ with coefficient $a_t = -1$ has a neighborhood $U$ with the same structure as an open toric surface. A divisor $\Delta\in \set{\Delta}$ has the property that $\Delta|_{U}$ extends to a trivial $\Q$-Cartier divisor $\overline{\Delta}$ on a toric compactification $\overline{U}$. Loosely speaking, the contribution to the $E$-function or elliptic genus of $(X,D+\varepsilon\Delta)$ coming from $U$ can be interpreted as the contribution from intersection data on $U$ coming from the index of a differential operator on $\overline{U}$ associated to the trivial divisor $\overline{\Delta}$. Moreover, the triviality of $\overline{\Delta}$ implies that this differential operator must vanish identically. These constraints are sufficient to prove that such intersection data on $U$ is independent of the choice of $\Delta$ in the limit as $\varepsilon\to 0$. 

\subsection{Outline}
The outline of this paper is as follows: In section (\ref{admiss}) we investigate a general class of normal surface pairs (called \emph{admissible pairs}) for which the stringy $E$-function, etc., are well-defined. Most importantly, we show that if $X$ is a normal surface with a $G$-action, and without strictly log-canonical singularities, then $(X,\Delta_{X/G})$ is admissible. This is of course a vital prerequisite to proving the McKay correspondence for normal surfaces in the non-log-terminal setting. In section (\ref{rigidity}) we review a useful rigidity theorem for an elliptic genus-like operator on toric varieties, which we will use later to show that our generalizations of the stringy $E$-function and elliptic genus are well-defined. Finally, in section (\ref{orb inv}) we define $E_{orb}$ and $Ell_{orb}$ for normal surfaces without strictly log-canonical singularities and prove the McKay correspondence for these invariants. We also discuss some closed formulas for some of these invariants and compute some examples.

\section{Admissible divisors on normal surfaces}\label{admiss}
\begin{defn}\label{admissible}\rm
Let $D = \sum_{I} a_i D_i$ be a smooth connected divisor with simple normal crossings on a smooth surface $X$. We say that $\sum_{I} a_i D_i$ is \emph{admissible} if it satisfies the following properties. For each $a_i = -1$: \\

$(1)$ $D_i \cong \Proj^1$.\\

$(2)$ Either $D_i$ intersects one component $D_{i_1}$ at a single point, or $D_i$ intersects two components $D_{i_1},D_{i_2}$ at a single point.
\end{defn}
The above definition is motivated by the following theorem in \cite{Veys}:

\begin{thm}\label{classification}
Let $p\in Y$ be normal surface singularity which is not log canonical. Let $\pi:X\rightarrow Y$ be the minimal log resolution of $p$. Then $\pi^{-1}(p) = \cup_{i\in I}D_i$ consists of the connected part $\cup_{i\in I,a_i<-1}D_i$ to which a finite number of chains are attached. If $D_0, D_1,...,D_r$ is such a chain, $D_0\subset K$, $D_\ell \cong \Proj^1$ for $1\leq \ell\leq r$, $a_1 \geq -1$ and $a_0< a_1 < \cdots < a_r < 0$. 
\end{thm}
Thus, an admissible divisor is a generalization of the exceptional divisor obtained by taking the minimal log resolution of a surface singularity which is not strictly log canonical. Note that the strictly log-canonical surface singularities constitute only a finite number of cases. See, for example, \cite{Aster}.

If $Y$ is a singular normal surface and $\Delta \subset Y$ a Weil divisor, we call $(Y,\Delta)$ an \emph{admissible pair} if there exists a log resolution $(X,D)\rightarrow (Y,\Delta)$ satisfying the following two properties: $(1)$ $D\subset X$ is admissible. $(2)$ For every component $D_i$ with $a_i = 0$, $D_i D_i\leq -2$. The reasoning behind this second criterion is to ensure that any two resolutions of $(Y,\Delta)$ satisfying $(1)$ and $(2)$ may be connected by a sequence of blow-ups and blow-downs $(X_i,D_i)$ so that each of the divisors $D_i$ in this sequence is admissible.

Admissibility appears to be a natural extension of the log-terminality condition for normal surfaces. For example, we have the following result, which will prove useful in the proof of the singular McKay correspondence.

Let $X$ be a normal surface without strictly log-canonical singularities, and let $G$ be a finite group acting on $X$. If $g:X\rightarrow X/G$ is the global quotient, then as in section (\ref{def McKay}), we let $\Delta_{X/G} \subset X/G$ be the $\Q$-Weil divisor satisfying $K_X = g^*(K_{X/G}-\Delta_{X/G})$. Then we have: 

\begin{prop}\label{McKay admissible}
$(X/G,\Delta_{X/G})$ is an admissible pair.
\end{prop}

We will require the following elementary fact from complex function theory:
\begin{lem}\label{automorphisms}
Let $H$ be a finite abelian group of automorphisms acting effectively on $\C$. Then $H$ is cyclic, given by rotation about a fixed point.
\end{lem}

\begin{proof}
Let $g\in H, g\neq e$. Then since $g$ lifts to an automorphism of $\Proj^1$ which fixes $\infty$, $g$ must have a single fixed point $p\in \C$. If $h\in H$ is any other element, then $gh\cdot p = hg\cdot p = h\cdot p$. Since $g$ has a unique fixed point in $\C$, it follows that $h\cdot p = p$. Thus, every element of $H$ has the same fixed point $p$. It follows that $H$ is a cyclic group given by rotations about $p$.
\end{proof}
We now proceed with the proof of proposition $(\ref{McKay admissible})$.
\begin{proof}
Let $(Y_0,E)$ be a minimal $G$-log resolution of $X$. If $E_0\subset E$ is a component with coefficient $b_0 = -1$, then by theorem (\ref{classification}), $E\cong \Proj^1$ and $E_0$ either intersects one component $E_1$ or two components $E_1, E_2\subset E$ at a single point. Note that if $g\in G$ fixes a point in $E_0$, then $gE_0 = E_0$, since the action of $G$ must permute divisor components with the same coefficients (and all divisor components with $-1$ coefficients are pair-wise disjoint). From this we can conclude that $gE_i = E_i$ if $g$ fixes $E_i\cap E_0, i=1,2$. Thus, any further blow-ups required to make $(Y,E)$ a $G$-normal resolution may be assumed to take place at points disjoint from $E_0$. We may therefore assume that $(Y_0,E)$ is already a $G$-normal pair.

By $G$-normality, the stabilizer of every point in $E$ is abelian. It follows that the subgroup $H\subset G$ that leaves $E_0$ invariant is abelian, since this subgroup must fix $E_0\cap E_1$. If $I_{E_0}\subset H$ is the cyclic inertia subgroup that fixes $E$, then $H/I_{E_0}=H_0$ is an abelian group acting effectively on $E_0$ and fixing $E_0\cap E_1$. If $H_0\neq \set{e}$, then by lemma (\ref{automorphisms}), $H_0$ is the cyclic group of rotations about $E_0\cap E_1$. In particular, either $H_0 =\set{e}$, in which case the image of $E_0$ under the map $Y_0\rightarrow Y_0/G$ is smooth, or $H_0$ fixes two points $p,q\in E_0$, and the image of $E_0\backslash\set{p,q}$ is smooth.

We now construct a resolution of singularities $(Z,\Delta)$ of $(X/G,\Delta_{X/G})$. After applying further blow-ups away from $E$, we may obtain a $G$-normal resolution of $(Y_0,E)$ with the property that every point in the resolution has an abelian stabilizer. For simplicity, we continue to refer to this resolution as $(Y_0,E)$. Since $Y_0$ has abelian stabilizers, the quotient variety $Y_0/G$ has at worst toric singularities. Let $(Z,\Delta)\rightarrow (X/G,\Delta_{X/G})$ be the smooth log resolution obtained by a minimal toric resolution of the singularities of $Y_0/G$.  After a finite sequence of equivariant blow-ups, we may obtain a resolution $(Y,D)\rightarrow (Y_0,E)$ and a smooth toriodal morphism $Y\rightarrow Z$ such that the following diagram commutes.
$$\begin{CD}
Y @>>> Z\\
@VVV       @VVV\\
X       @>>> X/G
\end{CD}$$
In the preceding analysis, we proved that for any component $E_0\subset Y_0$ with coefficient $a_0=-1$, we can find two points $p,q\in E_0$ so that the image of $E_0\backslash\set{p,q}$ under the quotient $Y_0\rightarrow Y_0/G$ is smooth. By our construction of $g:Y\rightarrow Y_0$, it follows that no curve in the exceptional set of $g$ gets mapped to $E_0\backslash\set{p,q}$. In particular, the proper transform of $E_0$ in $Y$ intersects at most two components of $D$; i.e., $(Y,D)$ remains an admissible pair.

We now prove that $(Z,\Delta)$ is admissible. Let $\cup D_j$ and $\cup \Delta_i$ denote the components of $D$ and $\Delta$. Let $b_j$ and $a_i$ denote their corresponding coefficients in $D$ and $\Delta$, repsectively. Since $f:(Y,D)\rightarrow (Z,\Delta)$ is a toroidal morphism with respect to the intersection data of $D$ and $\Delta$, $K_Y+\sum D_j = f^*(K_Z+\sum\Delta_i)$. Since we also have that $K_Y-D = f^*(K_Z-\Delta)$, we get that
\begin{align*}
f^*(K_Z-\sum a_i\Delta_i) &= 
f^*(K_Z+\sum\Delta_i-\sum (a_i+1)\Delta_i)\\
&=K_Y+\sum D_j - \sum (a_i+1)f^*\Delta_i\\
&=K_Y-\sum b_j D_j
\end{align*}
From the last two equations we can conclude that if $D_i$ maps to $\Delta_i$ with ramification $r_i$, then 
\begin{align}\label{pullback formula}
b_i = -1+r_i(a_i+1).
\end{align} 
It follows that $b_i = -1$ iff $a_i=-1$. Since $D$ is already admissible, we need only check that if $a_i = 0$, $\Delta_i \Delta_i \leq -2$. From the above equation, we must have that $b_i = r_i-1$. By theorem $(\ref{classification})$, this can only happen if $b_i = 0$. In other words, $\Delta_i$ must be an exceptional component of the toric resolution $Z\rightarrow Y_0/G$. Since we chose this resolution to be minimal, we must have that $\Delta_i\Delta_i \leq -2$. This completes the proof that $(X/G,\Delta_{X/G})$ is an admissible pair.
\end{proof}

\section{Rigidity and vanishing theorems on toric varieties}\label{rigidity}
Let $X$ be a smooth toric variety with toric divisors $D_1,...,D_\ell$ and big torus $T$. Consider the following $T$-equivariant vectorbundle over $X$, which we denote by $\mathcal{E}ll(a_1,...,a_\ell)$:
\begin{align*}
\bigotimes_{i=1}^\ell
   &\bigotimes_{n=1}^{\infty}
   \Lambda_{-y^{a_{i}}q^{n-1}}\mathcal{O}(-D_{i}) 
   \otimes \Lambda_{-y^{-a_{i}}q^{n}}\mathcal{O}(D_{i})\otimes
   S_{q^{n}}\mathcal{O}(-D_{i}) \otimes S_{q^{n}}\mathcal{O}(D_{i}) 
\end{align*}
The following theorem is proven in \cite{EllipALE}. Very similar results are also proven by Hattori in \cite{Hat}.
\begin{thm}
Assume $a_1D_1+...+a_\ell D_\ell = 0$ as a Cartier divisor for some choice of integers $a_1,...,a_\ell$. Then the $T$-equivariant index of $\mathcal{E}ll(a_1,...,a_\ell)$ vanishes identically.
\end{thm}

In fact, the theorem continues to hold if we merely assume that $a_i \in \Q$ and that $a_1D_1+...+a_\ell D_\ell = 0$ as a $\Q$-Cartier divisor. This is because for any integer $m$ satisfying $ma_i \in \Z$, we must have that $\chi_T(\mathcal{E}ll(ma_1,...,ma_\ell)) = 0$. By the fixed point formula, it is easy to see that $f(m)=\chi_T(\mathcal{E}ll(ma_1,...,ma_\ell))$ is a meromorphic function in the variable $m$. Since $f(m) = 0$ for infinitely-many $m$, we must have that $f(m)$ vanishes for all $m$, and in particular, for $m=1$. 

\begin{cor}
Let $(X,D)$ denote a toric Calabi-Yau pair. That is, $D = \sum a_i D_i$ is a $\Q$-Cartier toric divisor such that $K_X-D = 0$. Then 
$$Ell(X,D;\vec{t},z,\tau) \equiv 0.$$
\end{cor}

\begin{proof}
As proven in \cite{EllipALE}, $Ell(X,D;\vec{t},z,\tau)$ corresponds to the equivariant index of $\mathcal{E}ll(a_1+1,...,a_\ell+1)$, up to a normalization factor. Since $K_X = -D_1-...-D_\ell$ on a toric variety (see \cite{Fulton}), the Calabi-Yau condition implies that $\sum (a_i+1)D_i = 0$. It follows that $Ell(X,D;\vec{t},z,\tau)$ must vanish identically.
\end{proof}

The above results imply the following interesting rigidity theorem for the orbifold elliptic genus of a toric variety:

\begin{cor}\label{orb rigid}
Let $(X,D)$ be a toric Calabi-Yau pair and $G\subset T$ a finite subgroup. Then the $T$-equivariant orbifold elliptic genus $Ell_{orb}(X,D,G;\vec{t},z,\tau)$ vanishes.
\end{cor}

\begin{proof}
After replacing $(X,D)$ by a finite sequence of toric blow-ups $f:(\Bl{X},\Bl{D})\rightarrow (X,D)$, one can construct a smooth toric morphism $\mu:(\Bl{X},\Bl{D})\rightarrow (Y,\Delta)$ birational to the quotient by $G$ with $\mu^*(K_Y-\Delta) = K_{\Bl{X}}-\Bl{D} = f^*(K_X-D)$. By the equivariant McKay correspondence \cite{RW_McKay} and functoriality of the equivariant orbifold elliptic genus with respect to blow-ups, the equivariant orbifold elliptic genus of $(X,D,G)$ equals the equivariant elliptic genus of $(Y,\Delta)$. The fact that $(X,D)$ is a toric Calabi-Yau pair implies that $(Y,\Delta)$ is also Calabi-Yau. It follows that $Ell_{orb}(X,D,G;\vec{t},z,\tau)=0$.
\end{proof}
An interesting immediate consequence of the above result is a vanishing theorem for the orbifold $E$-function of a toric Calabi-Yau pair (where we again assume that the orbifold is a quotient by a finite subgroup $G\subset T$). To see this, note first that the vanishing of the orbifold elliptic genus of a toric Calabi-Yau pair $(X,D)$ implies the vanishing of the orbifold $\chi_y$ genus of $(X,D)$. However, it is easy to verify that the orbifold $E$-function of a toric pair $(X,D)$ is obtained from the orbifold $\chi_y$ genus by setting $y = uv$. Summarizing:
\begin{cor}\label{E rigid}
Let $(X,D)$ be a toric Calabi-Yau pair, and $G$ a finite subgroup of the maximal torus. Then the orbifold $E$-function, $E_{orb}(X,D,G;u,v) = 0$.
\end{cor}
\section{Orbifold invariants and the McKay correspondence}\label{orb inv}
To define the singular orbifold elliptic genus or stringy orbifold $E$-function in the non-log-terminal setting, we first note that the formulas for the orbifold elliptic genus and $E$-function for any smooth $G$-normal pair $(X,D)$ make sense as long as all the coefficients $a_i$ of the components of $D$ are not equal to $-1$. In the more general case where $(X,D)$ is an admissible pair, it is therefore natural to attempt to define $Ell(X,D,G;z,\tau)$ and $E(X,D,G;u,v)$ by introducing a perturbation $a_i+\varepsilon b_i$ to the coefficients of $D$, and then declaring the orbifold elliptic genus and orbifold $E$-function of $(X,D)$ to be the limit as the perturbation parameter tends to zero. 

As we will show below, the admissibility criterion guarantees that such a limit always exists. In general though, the value of the limit will depend on the choice of the perturbation. To give an example on the level of the stringy $\chi_y$ genus, suppose that $D_0$ has coefficient $a_0 = -1$ and intersects two divisors $D_1,D_2$ with coefficients $a_1$ and $a_2$. Consider a perturbation $a_i+\varepsilon b_i$ of these coefficients. Then the contribution of $D_0$ to the stringy $\chi_y$ genus of $(X,\sum (a_i+\varepsilon b_i) D_i)$ is equal to:
\begin{align*}
&\chi_y(D_0\cap D_1)\frac{(y-1)^2}{(y^{\varepsilon b_0}-1)(y^{1+a_1+\varepsilon b_1}-1)}+\\
&\chi_y(D_0\cap D_2)\frac{(y-1)^2}{(y^{\varepsilon b_0}-1)(y^{1+a_2+\varepsilon b_2}-1)}+
\chi_y(D_0^o)\frac{y-1}{y^{\varepsilon b_0}-1}\\
&=\frac{(y-1)^2(y^{\varepsilon(b_1+b_2)}-1)}
{(y^{\varepsilon b_0}-1)(y^{1+a_1+\varepsilon b_1}-1)(y^{1+a_2+\varepsilon b_2}-1)}
\end{align*}
In the last equality, we have used the relation $a_1+a_2 + 2 = 0$ which follows from the adjunction formula. It is easy to see that the limit as $\varepsilon\to 0$ in the above formula depends on the choice of $b_1,b_2$, and $b_3$. Notice, however, that if $b_1+b_2 = mb_3$ for some positive integer $m$, then 
\begin{align*}
&\frac{(y-1)^2(y^{\varepsilon(b_1+b_2)}-1)}
{(y^{\varepsilon b_0}-1)(y^{1+a_1+\varepsilon b_1}-1)(y^{1+a_2+\varepsilon b_2}-1)}=\\
&\frac{(y-1)^2(y^{\varepsilon b_0(m-1)}+...+1)}
{(y^{1+a_1+\varepsilon b_1}-1)(y^{1+a_2+\varepsilon b_2}-1)}
\end{align*}
and the limit of this expression as $\varepsilon\to 0$ depends only on $m$. Thus, one might attempt to solve this perturbation problem by requiring the perturbation coefficients to satisfy $b_1+b_2 = mb_0$ for some appropriate choice of $m$. In fact, one can show that the only choice of $m$ which makes the corresponding stringy $\chi_y$ genus invariant under blow-up is $m = -D_0 D_0$. In this case, one recovers the stringy $\chi_y$ genus of Willem Veys \cite{Veys}. It is less clear, however, that this procedure continues to work for more exotic invariants, such as the equivariant orbifold elliptic genus. As we will see however, the feasibility of this approach for all stringy invariants may ultimately be explained by rigidity phenomena associated to the equivariant orbifold elliptic genus.

The choice $b_1+b_2 = -D_0 D_0 b_0$ is clearly equivalent to requiring $c_1(b_1D_1+b_2D_2 + b_0 E_0)|_{E_0} = 0$. In what follows, we will see that this is the crucial property of a perturbation divisor which makes the $\varepsilon\to 0$ limit of a perturbed stringy invariant well-defined and independent of the choice of perturbation.

\begin{defn}\label{null perturbation}
\rm Let $(X,\sum a_i D_i)$ be a smooth admissible pair. We say that $\Delta_{\varepsilon}$ is a null-perturbation if $\Delta_{\varepsilon} = \sum \varepsilon b_i D_i$ and for any $a_j=-1$, $b_j\neq 0$, and $c_1(\Delta_\varepsilon)|_{D_j} = 0$.
\end{defn}

We now proceed to define the orbifold elliptic genus of $G$-normal divisor pairs.
\begin{defn}\rm
Let $(X,D)$ be a $G$-normal admissible divisor pair, and let $\Delta_\varepsilon$ be a null-perturbation. Then we define
\begin{align*}
Ell_{orb}(X,D,G;z,\tau)=
\lim_{\varepsilon\to 0}Ell_{orb}(X,D+\Delta_{\varepsilon},G;z,\tau)
\end{align*}
Similarly, we define the orbifold $E$-function
\begin{align*}
E_{orb}(X,D,G;u,v) = \lim_{\varepsilon\to 0}
E_{orb}(X,D+\Delta_\varepsilon,G;u,v)
\end{align*}
\end{defn}

If $(X,D)$ admits a $T$-action which commutes with $G$ and leaves $D$ invariant, we define the equivariant orbifold elliptic genus $Ell_{orb}(X,D,G;\vec{t},z,\tau)$ similarly. In order to make sense of the above definitions, we need to verify that $(1)$ the limits exist, and $(2)$ that the limits are independent of the choice of null-perturbation. To do this, we will make use of the rigidity theorems given in section $(\ref{rigidity})$. To summarize, we wish to prove the following proposition:

\begin{prop}\label{well-defined}
Let $(X,D)$ be $G$-normal, admissible. Assume that $X$ admits a torus action $T$ that commutes with the action of $G$ and preserves $D$. Let $\Delta_\varepsilon$ be a null-perturbation. Then
\begin{align*}
\lim_{\varepsilon\to 0}Ell_{orb}(X,D+\Delta_{\varepsilon},G;\vec{t},z,\tau)
\end{align*}
and
\begin{align*}
\lim_{\varepsilon\to 0}
E_{orb}(X,D+\Delta_\varepsilon,G;u,v)
\end{align*}
exist and are independent of the choice of null-perturbation.
\end{prop}

\begin{proof}
For simplicity, we may assume that $T=S^1$. Note that we allow for the possibility that the action is trivial. Let $\Delta_\varepsilon\subset X$ be a null-perturbation.

Fix a divisor component $D_t$ with $a_t=-1$. Let $H\subset G$ be the subgroup that leaves $D_t$ invariant. As shown in the proof of proposition $(\ref{McKay admissible})$, $H$ is abelian. Now $D_t$ has a $T\times H$-invariant open neighborhood isomorphic to the total space of the bundle $\mathcal{O}(-m_t)\rightarrow \Proj^1$. If we regard this neighborhood as a toric variety, it is easy to see that $T$ and $H$ act on $\mathcal{O}(-m_t)$ as subgroups of the maximal torus. We therefore reduce our analysis to the toric situation as follows:

Let $D_{t_k}, k\leq 2$ be the divisors intersecting $D_t$ with coefficients $a_{t_k}$. Consider the toric variety $\Proj^1\times\Proj^1$ with toric divisors $D_1 = \set{0}\times \Proj^1$, $D_2 = \Proj^1\times \set{0}$, $D_3 = \set{\infty}\times\Proj^1$, and $D_4=\Proj^1\times\set{\infty}$. Assign coeffients $a_{t_1}, a_{t_2}$ to $D_1$ and $D_2$, and coefficients $-a_{t_1}-2, -a_{t_2}-2$ to $D_3, D_4$. Then all these coefficients are distinct from $-1$, and 
\begin{align*}
&K_{\Proj^1\times\Proj^1}-a_{t_1}D_1-a_{t_2}D_2+(a_{t_1}+2)D_3+(a_{t_2}+2)D_4\\ 
&\equiv K_{\Proj^1\times\Proj^1}-D_{\Proj^1\times\Proj^1} = 0
\end{align*}
Blowing up $\Proj^1\times\Proj^1$ at $D_1\cap D_2$ produces an exceptional divisor $E\cong \Proj^1$ with discrepancy coefficient equal to $-1$. Let $p\in E$ be a fixed point of the action by the maximal torus. Let $Y$ be the toric variety obtained by blowing up at $p$ $m_t-1$ times. Then the proper transform of $E$ (which we continue to refer to as $E$) has discrepancy coefficient equal to $-1$ with regard to the map $f:Y\rightarrow \Proj^1\times\Proj^1$. Define $D_Y$ to be the toric divisor satisfying $K_Y-D_Y = f^*(K_{\Proj^1\times\Proj^1}-D_{\Proj^1\times\Proj^1})$. It is easy to see that $E$ intersects two components $E_1,E_2$ of $D_Y$ with coefficients $a_{t_1}$ and $a_{t_2}$. Moreoever, as for $D_t$, $E$ has a toric open neighborhood isomorphic to the total space of the bundle $\mathcal{O}(-m_t)\rightarrow\Proj^1$. Since both $T$ and $H$ may be regarded as subgroups of the maximal torus of $\mathcal{O}(-m_t)$, the $T$ and $H$ actions on this neighborhood extend naturally to all of $Y$. 

Let $\Delta_{\varepsilon}|_{D_t} = \varepsilon( b_1D_{t_1}+b_2 D_{t_2}+b D_t)|_{D_t}$. The condition that $c_1(\Delta_{\varepsilon})|_{D_t}=0$ implies that the divisor $b_1E_1+b_2 E_2 + bE$ is linearly equivalent to zero over the open neighborhood $\mathcal{O}(-m_t)$ of $E$. The set of toric Cartier divisors linearly equivalent to zero is in one-to-one correspondence with the set of linear functionals in the dual lattice of the toric variety. It follows that $b_1E_1+b_2 E_2+bE$ is the restriction to $\mathcal{O}(-m_t)$ of a Cartier divisor on all of $Y$ that is linearly equivalent to zero. Denote this divisor by $\Delta_Y$.

It is easy to see that the contribution coming from $D_t$ to the expression for $Ell_{orb}(X,D+\Delta_{\varepsilon},G;\vec{t},z,\tau)$ is equal to $[G:H]$ times the contribution coming from $E$ to the expression for $Ell_{orb}(Y,D_Y+\varepsilon\Delta_Y,H;\vec{t},z,\tau)$. Denote this contribution by $F(\varepsilon,\vec{t},z,\tau)$, and the remaining terms in the expression for $Ell_{orb}(Y,D_Y+\varepsilon\Delta_Y,H;\vec{t},z,\tau)$ by $G(\varepsilon,\vec{t},z,\tau)$. Since $K_Y-D_Y-\varepsilon\Delta_Y = 0$, corollary $(\ref{orb rigid})$ implies that
\begin{align*}
&Ell_{orb}(Y,D_Y+\varepsilon\Delta_Y,H;\vec{t},z,\tau)\\
&= F(\varepsilon,\vec{t},z,\tau)+G(\varepsilon,\vec{t},z,\tau)\\
&= 0
\end{align*}
for all $\varepsilon$. Since $G(\varepsilon,\vec{t},z,\tau)$ does not involve any divisor terms with $-1$ coefficients, $\lim_{\varepsilon\to 0}G(\varepsilon,\vec{t},z,\tau)$ exists and is independent of the coefficients of $\Delta_Y$. It follows that $\lim_{\varepsilon\to 0}F(\varepsilon,\vec{t},z,\tau)$ exists and is independent of the coefficients of $\Delta_Y$. This completes the proof for the case of the orbifold elliptic genus. The case for the orbifold $E$-function follows the same analysis upon applying corollary $(\ref{E rigid})$ in place of corollary $(\ref{orb rigid})$.
\end{proof}

As in the log-terminal setting, the orbifold elliptic genus and stringy $E$-function of an admissible pair satisfies the following functoriality property with respect to birational morphisms:
\begin{thm}
Let $(Y,D_Y)$ be a smooth, $G$-normal, $T$-equivariant admissible pair and $f:(X,D_X)\rightarrow (Y,D_Y)$ a $G\times T$ equivariant blow-up. Then 
\begin{align*}
&Ell_{orb}(X,D_X,G;\vec{t},z,\tau)=
Ell_{orb}(Y,D_Y,G;\vec{t},z,\tau)\\
&E_{orb}(X,D_X,G;u,v)=
E_{orb}(Y,D_Y,G;u,v)
\end{align*}
\end{thm}
\begin{proof}
We prove the case for the $E$-function. The case for the orbifold elliptic genus is exactly the same. Let $\Delta_{\varepsilon}\subset Y$ be a null-perturbation with respect to the divisor $D_Y$. Since $f^*\Delta_\varepsilon\cdot D_t = \Delta_{\varepsilon}\cdot f_* D_t$, $f^*\Delta_\varepsilon$ is a null-perturbation with respect to $D_X$. Since $K_X-D_X-f^*\Delta_{\varepsilon} = f^*(K_Y-D_Y-\Delta_\varepsilon)$,
\begin{align*}
E_{orb}(X,D_X+f^*\Delta_\varepsilon,G;u,v)=E_{orb}(Y,D_Y+\Delta_\varepsilon,G;u,v)
\end{align*}
by functoriality of the orbifold $E$-function for pairs. Taking the limit as $\varepsilon\to 0$ completes the proof.
\end{proof}

We may consequently define the \emph{singular orbifold elliptic genus } and \emph{singular orbifold $E$-function} for any $G$-equivariant admissible pair $(Y,\Delta)$ as follows: For $(X,D)\rightarrow (Y,\Delta)$ a $G$-normal admissible resolution, we define 
$$Ell_{orb}(Y,\Delta,G;z,\tau) = Ell_{orb}(X,D,G;z,\tau).$$
We define the orbifold $E$-function and equivariant orbifold elliptic genus of $(Y,\Delta)$ similarly. By the above proposition, the definition is independent of the choice of admissible resolution. 

\subsection{Singular McKay correspondence}
First, recall the following definition of a toroidal morphism:
\begin{defn}\rm
Let $\mu:(Y,D)\rightarrow (Z,\Delta)$ be a map of smooth varieties with simple normal crossing divisors. Note that the intersection data of the components of $D$ and $\Delta$ induce a stratification on $Y$ and $Z$.
We say that $\mu$ is a \emph{toroidal morphism} if:\\

\noindent $(1)$: $\mu : Y\backslash D \rightarrow Z\backslash\Delta$ is an unramified cover.\\

\noindent $(2)$: $\mu$ maps the closure of a stratum in $Y$ to the closure of a 
stratum in $Z$.\\

\noindent $(3)$: Let $U_z$ be an analytic neighborhood of $z \in Z$ such that the 
components of $\Delta$ passing through $z$ correspond to coordinate 
hyperplanes. Then for $y \in \mu^{-1}(z)$, there exists an analytic 
neighborhood $U_u$ of $y$ such that the components of $D$ passing 
through $y$ correspond to coordinate hyperplanes of $U_y$. Moreover, 
the map $U_y \rightarrow U_x$ is given by monomial functions in the 
coordinates. 
\end{defn}

We are now in a position to state and prove the singular McKay correspondence for surfaces in the non-log-terminal setting.
\begin{thm}
Let $G$ be a finite group acting holomorphically on a normal surface $X$. Assume $X$ does not have strictly log-canonical singularities. Moreoever, let $T$ be a compact torus that acts on $X$ and commutes with the action of $G$. Then
\begin{align*}
Ell_{orb}(X,G;\vec{t},z,\tau) = Ell(X/G,\Delta_{X/G},\vec{t},z,\tau)
\end{align*}
and
\begin{align*}
E_{orb}(X,G;u,v) = E_{str}(X/G,\Delta_{X/G};u,v)
\end{align*}
\end{thm}
\begin{proof} 
Let $\psi: X\rightarrow X/G$ be the global quotient map. As in the proof of proposition $(\ref{McKay admissible})$, we may construct the following commutative diagram:
$$\begin{CD}
(Y,D) @>\mu >> (Z,\Delta)\\
@VVV       @VVV\\
(X,0)       @>\psi >> (X/G,\Delta_{X/G})
\end{CD}$$
Here the vertical maps are resolutions of singularities. It is evident from the construction of the above diagram in proposition $(\ref{McKay admissible})$ that the map $\mu:(Y,D)\rightarrow (Z,\Delta)$ is a toroidal morphism. Moreover, all the maps are equivariant with respect to the $T$-action. 

Let $\Delta_\varepsilon\subset Z$ be a null-perturbation with respect to $\Delta$. Then equation $(\ref{pullback formula})$ implies that $\mu^*\Delta_\varepsilon$ is a null-perturbation with respect to $D$. It suffices to prove the equations
\begin{align*}
&Ell_{orb}(Y,D+\mu^*\Delta_\varepsilon,G;\vec{t},z,\tau)=
Ell(Z,\Delta+\Delta_\varepsilon;\vec{t},z,\tau)\\
&E_{orb}(Y,D+\mu^*\Delta_\varepsilon,G;u,v)=
E_{str}(Z,\Delta+\Delta_\varepsilon;u,v)
\end{align*}
For then the singular McKay correspondence follows after taking the limit as $\varepsilon\to 0$ and applying functoriality of the elliptic genera and $E$-functions with respect to the vertical arrows in the above diagram. 

We have thus reduced the problem to proving
\begin{align}\label{toroidal}
&Ell_{orb}(Y,D',G;\vec{t},z,\tau)=
Ell(Z,\Delta';\vec{t},z,\tau)\\
&E_{orb}(Y,D',G;u,v)=
E_{str}(Z,\Delta';u,v)
\end{align}
where $\mu:(Y,D')\rightarrow (Z,\Delta')$, and where $D'$, $\Delta'$ are $\Q$-divisors with the same components as $D$ and $\Delta$, \emph{but with no $-1$ coefficients}. For notational simplicity, we continue to write $D$ and $\Delta$ for $D'$ and $\Delta'$.

The case for the equivariant orbifold elliptic genus follows the exact same argument as the proof of theorem $(7)$ in \cite{RW_McKay}.
We therefore prove the case for the $E$-function. Let $D_j$, $j\in J$ denote the components of $D$, and $\Delta_i$,  $i\in I$ denote the components of $\Delta$. Let $Y^o = Y\backslash\cup_j D_j$, and define $Z^o$ similarly. By our construction of $(Y,D)$ and $(Z,\Delta)$, if $Y^g\neq \emptyset$ for any $g\neq e$, then $Y^g$ must intersect a component of $D$. It follows that the contribution to $E_{orb}(Y,D,G;u,v)$ coming from $\emptyset\subset J$ is given by $E(Y^o/G;u,v)$. Since the map $\mu:Y^o\rightarrow Z^o$ is simply the quotient by the free action of $G$, we have that $E(Y^o/G;u,v) = E(Z^o;u,v)$. This is precisely the contribution to $E_{str}(Z,\Delta;u,v)$ coming from $\emptyset\subset I$.

Next consider the contributions to $E_{orb}$ and $E_{str}$ coming from one-element subsets of $J$ and $I$. The contribution to $E_{str}(Z,\Delta;u,v)$ coming from the subset $\set{\Delta_i}$ is clearly given by
\begin{align*}
E(\Delta_i^o;u,v)\frac{uv-1}{(uv)^{a_i+1}-1}
\end{align*}
where $a_i$ is the coefficient of $\Delta_i$ in $\Delta$. Let $D_i\subset D$ be a component that maps to $\Delta_i$ under $\mu$. Then the coefficient $b_i$ of $\Delta_i$ is $-1+r_i(a_i+1)$, where $r_i$ is the order of the cyclic inertia subgroup $\Lambda_i\subset G$ that stabilizes $D_i$. The contribution to $E_{orb}(Y,D,G;u,v)$ coming from $G$-orbits of $\set{D_i}$ and conjugacy classes $[g], g\in \Lambda_i$ is equal to:
\begin{align}\label{contribution}
\sum_{k=0}^{r_i-1}(uv)^{k(a_i+1)}E(D_i^o/C_k;u,v)\frac{uv-1}
{(uv)^{r_i(a_i+1)}-1}
\end{align}
Here $C_k$ is the subgroup of $G$ that leaves $D^i$ invariant and commutes with $k\in Z_{r_i}\cong \Lambda_i$. In fact, $C_k = \mathrm{Inv}_G(D_i)$ for all $k$. For if $k\in \Lambda_i$ and $x\in\mathrm{Inv}_G(D_i)$, then $x^{-1}hxh^{-1}$ clearly equals the identity map on $D_i$ and acts trivially on the normal bundle to $D_i$. It follows that $x^{-1}hxh^{-1} = e$ since the action of $G$ is effective. Since $D_i^o/\mathrm{Inv}_G(D_i)=\Delta_i^o$, $(\ref{contribution})$ reduces to:
\begin{align}\label{contrib}
\sum_{k=0}^{r_i-1}(uv)^{k(a_i+1)}E(\Delta_i^o;u,v)\frac{uv-1}
{(uv)^{r_i(a_i+1)}-1}
\end{align}
Making use of the identity 
$$(uv)^{r_i(a_i+1)}-1=((uv)^{a_i+1}-1)\bigg\{\sum_{k=0}^{r_i-1}(uv)^{k(a_i+1)}\bigg\}$$ 
we may therefore identify expression $(\ref{contrib})$ with $E(\Delta_i^o;u,v)\frac{uv-1}{(uv)^{a_i+1}-1}$.

It remains to show that the contribution to $E_{str}(Z,\Delta;u,v)$ coming from two-element subsets $I'\subset I$ is equal to the contribution to $E_{orb}(Y,D,G;u,v)$ coming from two element subsets $J'\subset J$, and from conjugacy classes $[g]$ that fix isolated points. 

Since $(\ref{toroidal})$ holds for the orbifold elliptic genus, it must hold also for the $\chi_y$ genus. Setting $v=1$ in the above computations proves also that the contribution to $\chi_y(Z,\Delta)$ coming from $\emptyset$ and one-element subsets of $I$ corresponds to the contribution to $\chi_y(Y,D,G)$ coming from $\emptyset$, and from  $G$-orbits of one-element subsets together with the the conjugacy classes of elements in the corresponding inertia subgroups. These two facts together imply that the contribution to $\chi_y(Z,\Delta)$ coming from two-element subsets is equal to the contribution to $\chi_y(Y,D,G)$ coming from two-element subsets and from conjugacy classes $[g]$ that fix isolated points. Since both such contributions involve summations over isolated points, by setting $y=uv$ we obtain the formula for the contribution of these points to the orbifold and stringy $E$-functions. This completes the proof of $(\ref{toroidal})$ for the $E$-function. 
\end{proof}

\subsection{Closed formulas for stringy invariants}
We compute a closed expression for $E_{orb}(X,D,G;u,v)$, where $(X,D)$ is $G$-normal, admissible. Let $D = \sum_I a_i D_i$. Let $T\subset I$ index the set of all divisors with $a_t=-1$. For such a divisor $D_t$, let $D_{t_k}, k\leq 2$ denote the divisors which intersect $D_t$. For simplicity assume $k=2$; the case $k=1$ is analogous. Fix a null-perturbation $\Delta_\varepsilon=\sum\varepsilon b_i D_i$. Define $G_t\subset G$ to be the subgroup that fixes $D_t$. We first compute the contribution to $E_{orb}(X,D,G;,u,v)$ coming from $D_t$ and $g\in G_t$. Throughout, it will be convenient to make the change of variable $w=uv$.

Note first that if +$J\subset I^g$, then $D_t\cap D_J^o=\emptyset$ unless $J\subset\set{t,t_1,t_2}$. Moreover, $D_t\cap D_J^o/C(g,D_t,J)\cong D_t\cap D_J^o$. Therefore, the contribution coming from $D_t$ and $g$ is given by:
\begin{align*}
&\lim_{\varepsilon\to 0}
\sum_{J\subset\set{t,t_1,t_2}}
w^{F(g,D+\Delta_{\varepsilon})}E(D_J^o;u,v)\prod_J\frac{w-1}{w^{1+a_j+\varepsilon b_j}-1}\\
&= w^{F(g,D)}\frac{m_t(w-1)^2}{(w^{a_{t_1}+1}-1)
(w^{a_{t_2}+1}-1)}
\end{align*}
Here $m_t = -D_t D_t$. Next, suppose that $g$ acts invariantly on $D_t$ without fixing it. Then, as shown in the proof of proposition (\ref{McKay admissible}), $g$ acts on $D_t$ via rotation around the points $p_k =D_t\cap D_{t_k}$. Let $\alpha \in \Q\cap [0,1)$ be the infinitesimal weight of the $g$ action on the tangent space to $D_t$ at $p_1$. Then $1-\alpha$ is the corresponding infinitesimal weight at $p_2$. Let $\gamma_1,\gamma_2\in \Q\cap [0,1)$ denote the infinitesimal $g$-weights on the normal bundle to $D_t$ at $p_1,p_2$. 
For ease of notation, let $a = a_{t_1}+1$, $b_1 = b_{t_1}$, and $b_2 = b_{t_2}$.
Then the contribution from $g$ and $D_t$ to $E_{orb}(X,D,G;,u,v)$ is given by
\begin{align*}
\lim_{\varepsilon\to 0}\frac{(w-1)^2}{w^{\varepsilon b_t}-1}\bigg\{
\frac{w^{\alpha(a+\varepsilon b_1)+\varepsilon\gamma_1 b_t}}
{w^{a+\varepsilon b_1}-1}+
\frac{w^{(1-\alpha)(-a+\varepsilon b_2)+\varepsilon\gamma_2 b_t}}
{w^{-a+\varepsilon b_2}-1}
\bigg\}
\end{align*}
By a straight-forward but tedious computation, this limit evaluates to
\begin{align*}
w^{\alpha a}\big\{\alpha m_t+\gamma_1-\gamma_2+
w^a((1-\alpha)m_t+\gamma_2-\gamma_1)\big\}
{(w-1)^2}
\end{align*}
As a check, we may verify that this expression treats $p_1$ and $p_2$ on equal footing. Call the above expression 
$H(t,g)$. 

Now for any $g\in G$ let $T_1(g)$ index the set of $C(g)$-orbits of divisors $[D_t]$ such that $g\in G_t$. Let $T_2(g)$ index the set of $C(g)$-orbits of divisors $[D_t]$ such that $g \in \mathrm{Inv}_G(D_t)\backslash G_t$. Then, putting the above calculations together, we obtain the following closed expression for $E_{orb}(X,D,G;u,v)$:
\begin{align*}
&E_{orb}(X\backslash\cup_T D_t,D+\sum_T D_t,G;u,v)+\\
&\sum_{[g]}\bigg\{\sum_{[t]\in T_1(g)}
w^{F(g,D)}\frac{m_t(w-1)^2}{(w^{a_{t_1}+1}-1)(w^{a_{t_2}+1}-1)}
+\sum_{[t]\in T_2(g)} H(t,g)\bigg\}
\end{align*}
Note that the summands on the second line do not depend on the choice of representatives $g\in [g]$ or $t\in [t]$. Also, setting $G=\set{e}$, we recover Willem Veys' expression for the stringy $E$-function of the admissible pair $(X,D)$.

The above computation provides a direct verification of proposition $(\ref{well-defined})$ for the case of the orbifold $E$-function. For more exotic invariants, such as the orbifold elliptic genus (or its equivariant analogues), a corresponding verification by direct computation appears substantially more difficult. However, we do have the following explicit formula for the ordinary elliptic genus of a smooth admissible pair. Letting $x_j$ denote the formal Chern roots of $TX$, and using the same notation above for the coefficients of $D$, we have:
\begin{prop}
The elliptic genus of a smooth admissible pair $(X,D)$ is given by:
\begin{align*}
&\int_X\prod_{TX}\frac{\twopi{x_j}\vartheta(\twopi{x_j}-z)}
{\vartheta(\twopi{x_j})}
\prod_{I\backslash T}\jacob{D_i}{(a_i+1)}
\times\\
&\prod_T\frac{\vartheta(\twopi{D_t}+2z)\vartheta(z)}
{\vartheta(\twopi{D_t}+z)\vartheta(2z)}
+\sum_T m_t\frac{\vartheta(a_{t_1}z)\vartheta((a_{t_1}+2)z)}{\vartheta((a_{t_1}+1)z)^2}
\end{align*}
\end{prop}
Note that for ease of notation we omit the dependence on $\tau$ in the above formula. For a proof, see \cite{RW_Stringy}.

\subsection{Examples}
We give some examples of the above results which shed some light on what types of data we can expect $e_{str}(X,D)$ to encode in the non-log-terminal setting. 
Let $G\subset GL(3,\C)$ be a finite subgroup that acts effectively on $\Proj^2 = \Proj(x:y:z)$ and preserves a smooth curve $C$ of degree $d > 3$. Let $f$ be the defining equation for $C$ and $V_f\subset \C^3$ the affine hypersurface cut out by $f$. Then the McKay correspondence gives the following simple formula for $e_{str}(V_f/G,\Delta_{V_f/G})$:
\begin{cor}\label{cor_eff}
$e_{str}(V_f/G,\Delta_{V_f/G}) = \frac{(|G|+1)e(C/G)}{3-d}-d$
\end{cor}
\begin{proof}
Let $g:Y\rightarrow V_f$ be the proper transform of $V_f$ obtained by blowing up $\C^3$ at the origin. Then since $f$ defines a smooth curve in $\Proj^2$, $Y$ is a resolution of singularities with exceptional set $C$ and $K_Y+(d-2)C = f^*K_{V_f}$. Let $S^1$ act on $\C^3$ via the diagonal action. Then $S^1$ commutes with $G$ and descends to an action on $V_f$ with a single fixed point. The induced $S^1$-action on $Y$ acts freely on $Y\backslash C$ and fixes $C$. Since $G$ also acts freely on $Y\backslash C$, the orbifold Euler number of $(Y,(2-d)C,G)$ reduces to: 
\begin{align*}
e_{orb}(Y,(2-d)C,G) &=
\sum_{[g]}e(C^g/C(g))\frac{1}{3-d}
\end{align*}
Since the action of $G$ on $\Proj^2$ is effective, for all $g\neq e$, the fixed locus of $g$ must be a proper linear subspace. Hence, $g$ fixes only finitely-many points of $C$, since $C$ is not a linear subspace. In particular
\begin{align*}
\sum_{[g]}e(C^g/C(g))\frac{1}{3-d} = \frac{e_{orb}(C,G)}{3-d}=\frac{e(C/G,\Delta_{C/G})}{3-d}
\end{align*}
Here, the last equality follows from the ordinary McKay correspondence. An easy computation gives that $e(C/G,\Delta_{C/G}) = e(C/G)+B$ where $B =\sum(\nu_i-1)p_i$ is the degree of the branch divisor in the Riemann-Hurwitz formula for the map $C\rightarrow C/G$. Thus, 
\begin{align*}
e(C/G,\Delta_{C/G}) &= e(C/G)-e(C)+|G|e(C/G)\\
&= (1+|G|)e(C/G)+2g-2\\
&= (1+|G|)e(C/G)+d(d-3)
\end{align*}
The formula in the corollary then follows from the McKay correspondence.
\end{proof}
\begin{rmk}\rm
In the above formula, it is interesting to examine the case $d=3$ and $|G|=1$. In this case $e(C/G) =0$. If it were possible to cancel the expression $\frac{(1+|G|)e(C/G)}{3-d}$, we would arrive at the candidate value $-3 = e_{str}(V_f)$, where $V_f$ now has an isolated elliptic singularity. In fact, this is the answer one obtains for $e_{str}(V_f)$ if one applies the perturbation by an ample divisor approach to defining $e_{str}(V_f)$ \cite{RW_Stringy}\cite{BL_Sing}. This suggests that there may be a way to make sense of the McKay correspondence for surfaces with elliptic singularities.
\end{rmk}

\begin{rmk}\rm
For $d=1$, the above formula gives $e_{str}(V_f/G,\Delta_{V_f/G}) = |G|$. This is in agreement with the classical McKay correspondence, since the conditions on $G$ in the above proof force $G$ to be abelian for the case $d=1$.
\end{rmk}

Again, let $f(x,y,z)$ be a homogeneous polynomial of degree $d\neq 3$ which defines a smooth curve in $\Proj^2$. Let $V_f\subset \C^3$ be the affine hypersurface cut out by $f$. Let $G=\Z_n$ act on $\C^3$ by the diagonal action. Since $f$ is homogeneous, this action clearly descends to $V_f$, and we have the following simple formula for $e_{str}(V_f/G,\Delta_{V_f/G})$:

\begin{cor}
$e_{str}(V_f/G,\Delta_{V_f/G}) = nd = d\cdot e_{str}(\C^3/G,\Delta_{\C^3/G})$. 
\end{cor}

\begin{proof}
Again, we resolve $V_f$ by blowing up at the origin, acquiring the exceptional curve $C$. As in the proof of corollary $(\ref{cor_eff})$, we have the following formula:

\begin{align*}
e_{orb}(V_f,G)=\sum_{[g]}e(C^g/C(g))\frac{1}{3-d}
\end{align*}

Since $C$ is fixed by $G$, $e(C^g/C(g)) = e(C) = 2-2g = d(3-d)$. Since $G=\Z_n$, the first equality follows. The second equality follows from the ordinary McKay correspondence for the log-terminal variety $\C^3/G$.
\end{proof}

In fact, we can do a bit better and compute the stringy $\chi_y$ genus, i.e., $E_{str}(V_f/G,\Delta_{V_f/G};y,1)$. We have $E_{orb}(V_f,G;y,1)$ is equal to:
\begin{align*}
&\sum_{g}E(C;y,1)\frac{y-1}{y^{3-d}-1}
=n\cdot\frac{y-1}{y^{3-d}-1}(1-g)(1+y)
\end{align*}
Thus, $E_{str}(V_f/G,\Delta_{V_f/G};y,1) = n\cdot\frac{y^2-1}{y^{3-d}-1}\frac{d(3-d)}{2}$.


\end{document}